\documentclass[12pt]{article}
\usepackage{amssymb}
\usepackage{amsmath}
\usepackage{theorem}

\newcommand{\bQ}{{\mathbb Q}}

\newcommand{\bP}{{\mathbb P}}
\newcommand{\bZ}{{\mathbb Z}}
\newcommand{\cA}{{\mathcal A}}

\newcommand{\cE}{{\mathcal E}}

\newcommand{\cS}{{\mathcal S}}
\newcommand{\cT}{{\mathcal T}}
\newcommand{\cO}{{\mathcal O}}

\newcommand{\ra}{{\rightarrow}}

\newcommand{\Pic}{{\mathrm{Pic}}}

\newcommand{\Sym}{{\mathrm{Sym}}}

\newcommand{\Alb}{{\mathrm{Alb}}}

\newcommand{\Jac}{{\mathrm{Jac}}}
\newcommand{\Hom}{{\mathrm{Hom}}}
\newcommand{\Def}{{\mathrm{Def}}}

{\theorembodyfont{\rm}
\newtheorem{dfn}{Definition}[section]

\newtheorem{exam}[dfn]{Example}
\newtheorem{rem}[dfn]{Remark}
\newtheorem{lem}[dfn]{Lemma}
}

\newtheorem{theo}[dfn]{Theorem}
\newtheorem{prop}[dfn]{Proposition}

\author{Brendan Hassett and Yuri Tschinkel}
\title{\Large Abelian fibrations and rational 
points on symmetric products}

\date{\today}

\begin{document}

\maketitle

\section{Introduction}

Let $X$ be an algebraic variety defined over a number field $K$
and $X(K)$ its set of $K$-rational points. 
We are interested in properties of $X(K)$ 
imposed by the global geometry of $X$. 
We say that rational points on $X$ are potentially dense
if there exists a finite field extension $L/K$ such 
that $X(L)$ is Zariski dense. 
It is expected - at least for surfaces - that if
there are no finite \'etale covers of $X$ dominating a variety 
of general type then rational points on $X$ are potentially dense. 
This expectation complements the conjectures of 
Bombieri, Lang and Vojta predicting that rational points
on varieties of general type are always 
contained in Zariski closed subsets. 
This dichotomy holds for curves: the nondensity
for curves of genus $\ge 2$ is a deep theorem 
of Faltings and the potential  
density for curves of genus 0 and 1 
is classical. 
 
\

In higher dimensions there are at present 
no general techniques to prove nondensity.
Of course, potential 
density holds for abelian and unirational varieties. 
Beyond this, density results rely on the classification and
explicit projective geometry of the classes of 
varieties under consideration. 
In dimension two potential 
density is unknown for K3 surfaces 
with finite automorphisms and without elliptic fibrations
(see \cite{bogomolov-tschi-2}). 
In dimension 3 potential density is unknown, for example, for
double covers $W_2\ra\bP^3$ ramified in a smooth
surface of degree 6, for general conic bundles, as well 
as for Calabi-Yau varieties 
(for density results see \cite{harris-tschi}, 
\cite{bogomolov-tschi-1}). 
 
\

In this paper  
we study density properties of 
rational points on symmetric products $X^{(n)}=X^n/{\mathbb S}_n$. 
If $C$ is a curve of genus $g$ and $n>2g-2$
the symmetric product admits a bundle structure over the Jacobian 
$\Jac(C)$, with fibers projective spaces $\bP^{n-g}$.  
We see that in this case rational points on  $C^{(n)}$ are
potentially dense.  
Contrary to the situation for curves, 
we are not guaranteed to find many rational points on 
sufficiently high symmetric products of arbitrary surfaces. 
In Section \ref{sect:symm} we show that 
if the Kodaira
dimension of a smooth surface 
$X$ is equal to $k$ then the Kodaira dimension of $X^{(n)}$ is 
is equal to $nk$. 
This leads us to expect the behavior of rational points on 
$X^{(n)}$ and $X$  to be quite similar.
At the same time we observe that 
symmetric products of K3 surfaces admit (at least birationally) 
abelian fibrations over projective spaces. 
In fact, even symmetric squares of
certain (nonelliptic) K3 surfaces have the structure
of abelian surface fibrations over $\bP^2$. 
This is the starting point for proofs of potential density
of rational points. 
   
\

Let us emphasize that if $X$ is a variety over a number field $K$  
then Zariski density of rational points
on $X$ defined over degree $n$ field extensions of $K$
is not equivalent to Zariski density of $K$-rational 
points  on $X^{(n)}$. Of course, the first condition is
weaker than the second. Furthermore, if rational points on $X$
are potentially dense then they are potential dense
on $X^{(n)}$ as well.

\

This paper is organized as follows.  In Section \ref{sect:symm}
we recall general properties of symmetric products
and Hilbert schemes of surfaces.  Section 
\ref{sect:abelian} sets up generalities
concerning abelian fibrations ${\cal A}\ra B$. 
Potential density for ${\cal A}$ follows 
once one can find a ``nondegenerate''
multisection for which potential density holds.
In Section \ref{sect:elliptic} we prove widely-known results
concerning the existence of elliptic curves on K3 surfaces.  
Then we turn to potential density for symmetric products of K3 surfaces.
First, in Sections \ref{sect:density} and \ref{sect:sn}, we prove
potential density for sufficiently high symmetric powers
of arbitrary K3 surfaces.  
This is followed in Section \ref{sect:s2} with 
more precise results for symmetric squares of K3 surfaces of 
degree $2m^2$.

\
Throughout this paper, {\em generic} means `in a nonempty
Zariski open subset' whereas {\em general} means `in a
nonempty analytic open subset.'

\

{\bf Acknowledgements}. The first author was partially supported by
an NSF postdoctoral fellowship. 
The second author was partially supported by the NSA. 
We are very grateful to Joe Harris and Barry Mazur for
their help and encouragement.

\section{Generalities on symmetric products}
\label{sect:symm}

Let $X$ be a smooth projective variety over a field $K$.
Denote by $X^{n}=X\times_K ... \times_K X$ the $n$-fold
product of $X$. The symmetric group ${\mathbb S}_n$ 
acts on $X^n$. The quotient $X^{(n)}=X^n/{\mathbb S}_n$ is
a projective variety, called the symmetric product.

If $X$ has dimension 1 then $X^{(n)}$ is smooth and for
$n>2g-2$ the symmetric product  
$X^{(n)}$ is a projective bundle over the 
Jacobian $\Jac(X)$, 
with fibers projective spaces of dimension $n-g$ 
(see \cite{milne}, Ch. 4).   
In particular, rational points on $X^{(n)}$ are potentially
dense for $n>2g-2$. 

If $X$ has dimension 2 then $X^{(n)}$ is no longer smooth;
it has Gorenstein singularities (since the group action
factors through the special linear group).  
The Hilbert scheme of length $n$ zero-dimensional
subschemes is a crepant resolution of $X^{(n)}$ 
$$
\varphi\,:\, X^{[n]}\ra X^{(n)}
$$ 
(see \cite{beauville-83}, Section 6 and the references therein).
In particular, $\varphi^* \omega_{X^{(n)}} = \omega_{X^{[n]}}$.
The same holds for pluricanonical differentials. 
On the other hand, we have the isomorphism 
$$
H^0(X^n, \omega_{X^n}^m)^{{\mathbb S}_n}=H^0(X^{(n)},
\omega_{X^{(n)}}^m).
$$
We are using the fact that the quotient map 
$X^n\ra X^{(n)}$ is unramified away from a codimension two
subset and pluricanonical differentials extend over
codimension two subsets.
We conclude that pluricanonical differentials on the Hilbert
scheme correspond to ${\mathbb S}_n$-invariant 
differentials on the $n$-fold
product $X^n$:
$$
H^0(X^{[n]}, \omega_{X^{[n]}}^m)\simeq 
H^0(X^n, \omega_{X^n}^m)^{{\mathbb S}_n}.  
$$
Since
$$
H^0(X^n,\omega_{X^n}^m)^{{\mathbb S}_n}\simeq \Sym^n H^0(X,\omega_X^m)
$$  
we obtain the following:

\begin{prop}
Let $X$ be a smooth surface. If $X$ has Kodaira dimension $k$ then
$X^{(n)}$ has Kodaira dimension $nk$. 
\end{prop}

If $X$ is a K3 surface we can be more precise:  
$X^{[n]}$ is a holomorphic symplectic
manifold (see \cite{beauville-83}, Section 6).
In particular, the canonical bundle of $X^{[n]}$ remains
trivial. 

An important ingredient in the proofs of potential density is
the construction of a multisection of the abelian fibration.
The following proposition will help us verify that
certain subvarieties are multisections:

\begin{prop}\label{prop:*}
Let $X$ be a smooth projective surface and $C_1,..., C_n$ 
distinct irreducible curves. Consider the image $Z$ of 
$C_1\times ... \times C_n$ under the quotient map $X^n\ra X^{(n)}$.
The scheme-theoretic preimage $\varphi^{-1}(Z)\subset X^{[n]}$
has a unique irreducible component of dimension $\ge n$, denoted
by $C_1 * ... * C_n$. In particular, the homology class 
of $C_1 * ... * C_n$ is uniquely determined by the homology classes
of $C_1,..., C_n$.    
\end{prop}

{\em Proof.}
Let $(a_1,.., a_k)$ 
be a partition of $n$ and let ${\cal D}_{a_1,...,a_k}$
be the corresponding stratum in $X^{(n)}$. 
In particular, the  ${\cal D}_{a_1,...,a_k}$ are disjoint. 
The intersection of $Z$ with 
${\cal D}_{a_1,...,a_k}$ has dimension $\#\{ a_j\,|\, a_j = 1\}$. 
Each fiber of $\varphi$ over ${\cal D}_{a_1,...,a_k}$
is irreducible of 
dimension $\sum_{j=1}^k (a_j-1)$ (see \cite{briancon}).
It follows that the preimage of $Z\cap {\cal D}_{a_1,...,a_k}$
has dimension at most 
$$
\#\{ a_j\,|\, a_j = 1 \}  + \sum_{j=1}^k (a_j-1),
$$
which is less  than $n$, provided the $a_j$ are not all  equal to 1. 
$\square$

\section{Generalities on abelian fibrations}
\label{sect:abelian}

Let $\cA$ be an abelian variety defined over a field $K$
(not necessarily a number field).  
A point $\sigma \in \cA(K)$ is {\em nondegenerate} if 
the subgroup generated by $\sigma$ is Zariski dense in
$\cA$. 

\begin{prop}\label{prop:nondegenerate}
Let $\cA$ be an abelian variety over a number field $K$.
Then there exists a finite field extension $L/K$
such that $\cA(L)$ contains a nondegenerate point.
\end{prop}

{\em Proof.} 
We include an argument for completeness, since we could not find
a reference.  

\begin{lem} \label{lemm:rank}
Let $\cA$ be an abelian variety of dimension $\dim(\cA)$ defined
over a number field $K$. Then there exists a finite field extension $L/K$
such that the rank of the Mordell-Weil group $\cA(L)$ 
is strictly bigger than the rank of $\cA(K)$. 
\end{lem}

{\em Proof.} 
We first assume $\dim(\cA)>1$.  
Let $\Gamma$ be the {\em saturation} of $\cA(K)$ in $\cA({\bar K})$
(where ${\bar K}$ is the algebraic closure of $K$). 
This means that $\Gamma$ 
consists of all points $p$ such that a positive multiple of $p$ lies in
$\cA(K)$; in particular it contains all torsion points.
Find a smooth curve $C$ of genus $\ge 2$ 
in $\cA$, defined over a number field $K_1$.  By
Raynaud's version of the Manin-Mumford conjecture  
(see \cite{lang-3}, I 6.4 or \cite{oesterle} Theorem 1) 
we have that $C\cap \Gamma$ is finite.  
There exists a $L/K_1$ such that $C(L)$ contains 
a point $q$ outside $C\cap \Gamma$. 
It follows that $\cA(L)$ has higher rank. 
(This argument was communicated to us by B. Mazur.)

We now do the case of an elliptic curve $\cE$.  
Write $\cA=\cE \times \cE$ with projections $\pi_1$ and $\pi_2$;
we have $\cA(K)=\cE(K) \times \cE(K)$.  The argument above gives
a point $q\in \cA(L)$ not contained in the saturation of $\cA(K)$.
It follows that either $\pi_1(q)$ or $\pi_2(q)$ is not 
contained in the saturation of $\cE(K)$.  $\square$

\

We prove the proposition.  We may replace $\cA$ with an isogenous abelian
variety, so we may assume that $\cA$ is a product of geometrically
simple abelian varieties.
Our proof proceeds by induction on the number of simple components.  
Any nontorsion point $p$ of a geometrically simple abelian variety is
nondegenerate.  Indeed, Faltings' theorem implies that the
Zariski closure of ${\bZ}p$ is a finite union of abelian subvarieties.  
Hence it suffices to prove the inductive step:

\begin{lem}
Let $\cA_1$ and $\cA_2$ be abelian varieties over a number field $K$. 
Assume that $\cA_2$ is geometrically simple and 
$\cA_1$ and $\cA_2$ have nondegenerate $K$-points $p_1$ and $p_2$.
Then $\cA_1\times \cA_2$ has a nondegenerate point over some finite
extension $L/K$.
\end{lem}

{\em Proof.} 
For any pair of abelian varieties $\cA_1,\cA_2$ the 
group of homomorphisms $\Hom(\cA_1,\cA_2)$
is finitely generated as a module over $\bZ$. 
After a finite extension, we may assume
these are all defined over $K$.
More generally, we can consider the group generated
by abelian subvarieties of $\cA_1\times \cA_2$, or equivalently,
by homomorphisms from $\cA_1$ to $\cA_2$ defined only up to isogeny.  
This group equals $\Hom^0(\cA_1,\cA_2)= \Hom(\cA_1,\cA_2)\otimes \bQ$
(see \cite{mumford}, p. 172-176).

Assume that $(p_1,p_2)$ is contained in 
a proper abelian subvariety ${\cal B}$. 
We regard ${\cal B}$ as an element $\beta\in \Hom^0(\cA_1,\cA_2)$; 
in particular, $\beta(p_1)=p_2$.  We choose a $\bZ$-basis $(Z_1,...,Z_k)$ 
for $\Hom(\cA_1,\cA_2)$.  There exist integers $b_1,...,b_k,$ and $d\ne 0$
such that $(b_1 Z_1+...+b_kZ_k)(p_1)=dp_2$ in the Mordell-Weil group.
Hence $p_2$ is contained in the saturation of the 
subgroup of $\cA_2(K)$ generated
by the images of $p_1$ under the $Z_i$.  Conversely, if $q$
is not contained in this subgroup then $(p_1,q)$ is nondegenerate.
Applying Lemma \ref{lemm:rank}, we obtain a finite field extension
$L/K$ and a point $q\in \cA_2(L)$ with the desired property.  $\square$

\

Let $\cT$ be an $\cA$-torsor, defined over a field $K$.  
This means there is a $K$-isomorphism
$\mathrm{Aut}_0(\cT)\simeq \cA,$
so that for any $p\in {\cT}(M)$ (where $M$ is an extension of $K$),
the translation action induces an isomorphism
${\cA}(M) \ra {\cT}(M).$
There is a 1-1 correspondence between $M$-valued points of ${\cal T}$ 
and isomorphisms $\cA(M) \simeq \cT(M)$.

Consider the Albanese ${\Alb}(\cT)$ (see, for example, \cite{lang}
II. 3).  It is an abelian variety defined over $K$, such that 
there is a morphism $\cT \times \cT\ra {\Alb}(\cT)$ corresponding to
$(t_1,t_2) \rightarrow t_1-t_2$.  
For each zero-cycle of $\cT$, defined over $K$ and of degree zero,
we obtain a point in ${\Alb}(\cT)(K)$.  
Note that $\cA$ is naturally isomorphic to $\Alb(\cT)$
over $K$.  Indeed, the action of ${\cA}$ on ${\cT}$
induces an action on $\Alb(\cT)$, and the orbit of
zero is isomorphic to both ${\cA}$ and $\Alb(\cT)$.

\

Let $\pi\,:\, {\cal T}\ra B$ 
be an abelian fibration, that
is: ${\cal T}$ and  $B$ are normal and $B$ is connected and 
the fiber ${\cal T}_b$ 
over the generic point $b$ is a torsor for
an abelian variety $\cA_b$ over $K(b)$.
A multisection ${\cal M}$ of $\pi$ is the
closure of an $M$-valued point of ${\cal T}_b$, where $M$ is 
a finite field extension of $K(b)$ of degree $\deg(M)$. 
Let $\mu$ denote the generic point of ${\cal M}$.

\begin{prop}\label{prop:density}
Let ${\cal T}$ be an abelian fibration with  
multisections ${\cal M}_1$ and ${\cal M}_2$, 
both defined over a number field $K$. 
Let $\sigma$ be the distinguished section of 
${\cal T}(M_1)$. Denote by ${\cal M}_2'={\cal M}_2\times_B \mu_1$; it is 
a cycle defined over $M_1$. 
Assume that $K$-rational points on ${\cal M}_1$
are Zariski dense and that the cycle
$$\deg(M_2)\sigma - {\cal M}'_2$$
yields a nondegenerate point of $\Alb({\cal T})(M_1)$.   
Then $K$-rational points on ${\cal T}$ are Zariski dense.  
\end{prop} 

{\em Proof.} We restrict to an open subset of 
$B$ where $\pi $ and ${\cal M}$
are flat over $B$.  
Denote by ${\cal T}'$ the base change of ${\cal T}$ to 
${\cal M}_1$ with distinguished section $\sigma$.  
Consider the action of
${\cal A}(M_1)\simeq {\Alb}({\cal T})(M_1)$ on ${\cal T}(M_1)$.  
The translates
of $\sigma$ by the nondegenerate point 
are Zariski dense in ${\cal T}(M_1)$.
Note that translation by the 
nondegenerate point is defined over the number field $K$.
Each of these translates is 
birational to $\sigma $ over $K$ and therefore
its $K$-rational points are Zariski dense. 
The union of the closures in ${\cal T}'$ 
of these translates is Zariski dense;
hence $K$-rational points in ${\cal T}'$ are 
Zariski dense.  It remains to observe that ${\cal T}'$ dominates
${\cal T}$. $\square$

\begin{rem}
Our argument does not show that rational points are Zariski dense
in {\em any} fiber of ${\cal T}_p$, 
where $p$ is an $K$-rational point
of $B$. However, one knows that 
when the fibers are of dimension 1 then 
for $p'\in U_1(K)$, where $U_1$ is some
nonempty open subset of ${\cal M}_1$, the fibers over $p'$ 
have infinitely many $K$-rational points (see \cite{silverman}). 
Moreover, by a result of N\'eron, 
the rank of the Mordell-Weil group of
special fibers of abelian fibrations does not drop outside a thin
subset of points on the base of the fibration \cite{serre}. 
\end{rem}

\section{Elliptic families on K3 surfaces}
\label{sect:elliptic}

Throughout this section, we work 
over an algebraically closed field of characteristic 0.
An elliptic fibration is an abelian fibration of relative
dimension one. In the sequel an 
elliptic fibration dominating a K3 surface will be called an 
elliptic family.

The following theorem is attributed to Bogomolov and Mumford
(see \cite{mori-mukai}). We include a detailed proof
because it is crucial for our applications.

\begin{theo}\label{theo:fibrations}
Let $S$ be K3 surface and $f$ a divisor
class on $S$ such that $h^0({\cO}_S(f))>1$.  
Then there exists a smooth curve $B$ and 
an elliptic fibration ${\cal E}\ra B$
with the following properties:
\begin{enumerate}
\item   ${\cal E}$ dominates $S$;
\item   the generic fiber ${\cal E}_b$ is mapped
birationally onto its image;
\item   the class $f-{\cal E}_b$ is effective.
\end{enumerate}
\end{theo}

{\em Proof.} 
A genus one curve $C\subset S$ 
is a curve whose normalization $\tilde{C}$ is a
connected curve of genus one.  
It suffices to prove the result for a
singular curve $B$; we can always 
pull back to the normalization $\tilde{B}$. 

We may restrict to the case where $S$ 
is not an elliptic K3 surface.
We assume that $|f|$
has no fixed components (and thus no base points).
Indeed, if this is not the case then we extract the moving
part of $f$.  Since $S$ is not elliptic, we have $f^2>0$.  
We may also assume that the class $f$ is primitive;
otherwise, take the primitive effective generator $f'$ of ${\bZ}f$.
We still have $h^0({\cO}(f'))>1$ and $|f'|$ basepoint free
(again, using the fact that $S$ is not elliptic.)
See \cite{saint-donat} for basic results concerning linear series 
on K3 surfaces.

We shall use the following lemma, essentially proved in
\cite{mori-mukai}:
\begin{lem}
\label{Mori-Mukai}
For each $n>0$, a generic polarized K3 surface $(S_1,f)$
of degree $2n$ contains a one-parameter family of irreducible
curves with class $f$, such that the generic member is nodal 
of genus one.  
\end{lem}

{\em Proof.}
We first claim there exists a K3 surface $S_0$ containing two smooth rational
curves $D_1$ and $D_2$ meeting transversally at $n+2$ points.  
Let $S_0$ be the Kummer surface associated to the product of elliptic
curves $E_1$ and $E_2$, such that there exists an isogeny $E_1 \ra E_2$ of
degree $2n+5$.  Let $\Gamma$ be the graph of this isogeny and $p\in E_2$
a $2$-torsion point.  Now $\Gamma$ intersects 
$E_1 \times p$ transversally in $2n+5$ points, one of which 
is $2$-torsion in $E_1 \times E_2$.  
We take $D_1$ to be the image of $\Gamma$ and
$D_2$ to be the image of $E_1 \times p$;  $D_1$ and $D_2$ are smooth,
rational, and intersect transversally in $n+2$ points.  
The line bundle ${\cO}(f):={\cO}_{S_0}(D_1+D_2)$ is big and nef and thus 
has no higher cohomology (by Kawamata-Viehweg vanishing).

Let $\Delta$ be the spectrum of a discrete valuation
ring with closed point $0$ and generic point $\eta$.  
Let $\cS \rightarrow \Delta$ be a deformation of $S_0$ such that
$f$ remains algebraic.  We assume further that the class $f$ is 
ample and indecomposible in the monoid of effective curves in a
(geometric) generic fiber $S_1$.  These conditions are satisfied away 
from a finite union of irreducible divisors.  
Since $f$ has no higher cohomology, $D_1 \cup D_2$ is a specialization
of curves in the generic fiber and the deformation space 
$\Def(D_1 \cup D_2)$ is smooth of dimension $n+2$.  
Consider the locus in $\Def(D_1 \cup D_2)$ 
parametrizing curves with at least $\nu$ nodes;
this has dimension $\ge n+2-\nu$.  When $\nu=n+1$ the corresponding curves
are necessarily rational.  Each fiber of 
$\cS \rightarrow \Delta$ is not uniruled, and thus contains
a finite number of these curves.  In each fiber, the rational curves
with $n+1$ nodes deform to positive-dimensional families of curves
with $n$ nodes.  Hence $S_1$ contains a family of nodal
curves of genus one with the desired properties. $\square$

To complete the proof, we use a proposition suggested
by Joe Harris:

\begin{prop}
Let ${\cal S}\ra D$ be a projective morphism. 
Then there exists a scheme ${\cal K}_g({\cal S}/D)$ such that each
connected component is projective over $D$ and the fiber 
over each $d\in D$ is isomorphic
to the corresponding moduli space of 
stable maps ${\cal K}_g({\cal S}_d)$.
\end{prop}
 
{\em Proof.}
We refer to Kontsevich's moduli space of stable maps constructed in 
\cite{kontsevich-manin},\cite{fulton-rahul}. We first consider the 
special case when ${\cal S}=\bP^n_D$. Then 
$$
{\cal K}_g(\bP^n_D/D)= 
{\cal K}_g(\bP^n)\times D
$$
More generally, given an embedding  
${\cal S}\ra \bP^n_D$ over $D$, we define ${\cal K}_g({\cal S}/D)$ as
those elements of ${\cal K}_g(\bP^n_D/D)$ 
which factor through ${\cal S}$.
Since it is a closed subscheme it is projective over $D$.  $\square$

We finish the proof of Theorem \ref{theo:fibrations}.
There exists a projective family of K3 surfaces ${\cal S} \ra \Delta$
equpped with a divisor class $f$,
such that the (geometric) generic fiber satisfies the conditions of
Lemma \ref{Mori-Mukai} and the special fiber is $(S,f)$.  
Consider the component ${\cal K}_1({\cal S}/\Delta,f)$ of  
${\cal K}_1({\cal S}/\Delta)$ consisting of maps with 
image in the class $f$.  After a finite base change 
$\Delta' \ra \Delta$, there exists a geometrically irreducible
curve ${\cal C}_{\eta}\subset 
{\cal K}_1({\cal S}/\Delta,f)$ 
corresponding to an elliptic fibration dominating the generic
fiber ${\cal S}_{\eta}$.  
Let ${\cal C} \subset 
{\cal K}_1({\cal S}/\Delta,f)$ be the flat extension over $\Delta$
and ${\cal C}_0$ the corresponding flat limit.  

There may not be a `universal stable map' defined over
${\cal C}_0 \subset {\cal K}_1(S,f)$.  However, for each irreducible
reduced component $C_i \subset {\cal C}_0$, a universal stable map
exists after a finite cover $B_i \ra C_i$.  For some such $B_i$,
the resulting family of stable maps 
$\cE'_i \rightarrow B_i$ dominates $S$.
The image of the generic fiber contains a component of
genus one because no K3 surface is uniruled.  
$\square$

\section{Density of rational points}
\label{sect:density}

In this section $S$ denotes a K3 surface
defined over a number field $K$. 
Potential density holds for elliptic K3 surfaces and 
for all but finitely many families of K3 surfaces
with Picard group of rank $\ge 3$, and 
consequently for their symmetric products (see \cite{bogomolov-tschi-2}).
However, a 
general K3 surface has Picard group of rank 1. 
In the following sections
we will prove density results for symmetric
products of general K3 surfaces. 

\

By Theorem  \ref{theo:fibrations}, there is a family
of elliptic curves ${\cal E}$ dominating $S$. 
Let $E_1,\ldots,E_n$ be generic curves in the fibration
and assume that
$g=[E_i]$ is {\em big}; in particular, ${\cal E}$
is not an elliptic fibration on $S$. 
It follows that the general
member of $g$ is an irreducible curve of genus $>1$.
Note that we have a well defined class $g * ... * g$ in the homology 
of $S^{[n]}$, equal to 
the homology class of $C_1 * ... * C_n$, where the $C_i$ 
are irreducible curves in $g$ (see Proposition \ref{prop:*}).

\begin{theo}\label{theo:dense}
Let $S$ be a K3 surface satisfying the conditions 
of the previous paragraph.  
Assume that either
\begin{enumerate}
\item
$F=S^{[n]}$ admits
an abelian fibration ${\cal T}\ra B$ and $g * ... * g$ intersects the
proper transform of the generic
fiber positively, or
\item 
$F$ is birational to an abelian fibration, and $E_1* ... *E_n$
is a multisection.
\end{enumerate}
Then rational points on $F$ are potentially dense. 
\end{theo}

{\em Proof.} Throughout the proof, $L/K$ is some finite
field extension, which we will enlarge as necessary. 
We want to show that $L$-rational points are Zariski dense
on $F$. 

For generic smooth curves $C_1,..., C_n$ in $g$,
$C_1 * ... * C_n\subset S^{[n]}$ 
gives a multisection of ${\cal T}\ra B$. 
We denote this multisection by ${\cal M}_2$. 
This is clear under the first assumption.  Under
the second assumption, it follows from the fact that
$E_1 * ... * E_n$ is a multisection.

Choose a point $p\in B(L)$ corresponding to a 
smooth fiber ${\cal T}_p$ of ${\cal T}$.
Let ${\cal A}_p$ be the Albanese of ${\cal T}_p$. Choose a point
$m_1\in {\cal T}_p(L)$ so that the class of the cycle 
\begin{equation}\label{eqn:nondegenerate}
 \deg({\cal M}_2 )\cdot m_1 - {\cal M}_2|_{{\cal T}_p}
\end{equation}
is nondegenerate in ${\cal A}_p$ (see Proposition \ref{prop:nondegenerate}).
We may assume that $m_1$ corresponds to 
a subscheme $(s_1,...,s_n)\in S^{[n]}$ 
where the $s_i$ are distinct, 
$E_i$ contains $s_i$ for $i=1,...,n$, 
the $s_i$ and the $E_i$ are defined over $L$, 
and $L$-rational points are Zariski dense on each $E_i$. 
Then we have a multisection ${\cal M}_1$ for ${\cal T}$ 
given as (the proper transform of ) $E_1 * ... * E_n$. 
Note that $L$-rational points on ${\cal M}_1$ are Zariski dense.

It follows from Equation (\ref{eqn:nondegenerate}) 
that the pair $({\cal M}_1, {\cal M}_2)$
satisfies the 
nondegeneracy assumptions of Proposition \ref{prop:density}. 
Therefore, $L$-rational points are Zariski dense in $F$. $\square$ 

We employed two parallel sets of hypotheses because in some
applications the abelian fibration is only described over
the generic point of $B$, which makes intersection computations
difficult.  In other applications, the abelian fibration
is given by an explicit linear series, but the multisection
is difficult to control.   

\begin{rem}
Matsushita has proved a structure theorem for holomorphic symplectic
manifolds of dimension $2n$ admitting a 
fibration structure. In particular, he proved
that the base 
has dimension $n$, is Fano, 
has Picard group of rank 1, and log-terminal
singularities. Furthermore, the  
fibers admit finite \'etale covers which are
abelian varieties (see \cite{matsushita}). 
\end{rem}

\begin{rem}
We do not know how to produce abelian fibrations on 
symmetric products of Calabi-Yau varieties of dimension $\ge 3$.
For example, do they exist for quintic threefolds?
\end{rem}

\section{Potential density on $S^{[n]}$}
\label{sect:sn}

In this section we exhibit K3 surfaces $S$ defined
over a number field $K$ and  
satisfying the assumptions 
of Theorem \ref{theo:dense}.

\begin{theo}
Let $S$ be a K3 surface with
Picard group of rank 1 generated by  
a polarization $g$ of degree $2(n-1)$. 
Then 
there exists a finite extension $L/K$ such that
$L$-rational points
on $S^{[n]}$ are Zariski dense. 
\end{theo}

{\em Proof.}
Under our hypothesis, $g$ is basepoint free and yields
a morphism $S\ra \bP^{n}$ which is finite onto its
image. 
Furthermore, the generic member of $|g|$ is 
smooth of genus $n$ (see \cite{saint-donat}).

There is an abelian fibration 
over $B\subset \bP^n$, 
where $B$ corresponds to the locus of smooth curves in $|g|$.
Indeed, ${\cal T}\ra B$ is the degree $n$ component of the
relative Picard fibration
(see \cite{beauville-97}).
We claim that $S^{[n]}$ is
birational to ${\cal T}$. Given generic points $s_1,..., s_n$
on $S$ there is a smooth curve $C\in |g|$ passing through those
points. The line bundle ${\cal O}_{C}(s_1+...+s_n)$ 
is a generic point of $\Pic_n(C)$, and such a line bundle 
has a unique representation as an effective divisor.  
(We are using the fact that $C^{[n]}$ is birational
to $\Pic_n(C)$.)

To apply Theorem \ref{theo:dense} we must 
verify that (the proper transform of) $E_1 * ... * E_n$ 
is a multisection for ${\cal T}$. 
A generic curve $C\in |g|$  
intersects the 
union of the $E_i$ transversally in $n(2n-2)$ points. 
Under these assumptions, every subscheme parametrized
by
$C^{[n]}\cap (E_1 * ... * E_n)$ is reduced and there are
finitely many such subschemes. 
It particular, $C^{[n]}$ intersects 
$E_1 * ... * E_n$
in finitely many points. $\square$ 

\begin{rem}\label{rem:main}
The same argument applies if the rank of the Picard group
of $S$ is $>1$. Our proof uses only that 
$g$ is big and contains the class of an irreducible 
elliptic curve on $S$.  Under these conditions $|g|$
is basepoint free;  the base locus of any linear series
on a K3 surface has pure dimension one (see \cite{saint-donat}).
\end{rem}

\begin{theo}\label{theo:main}
Let $S$ be a K3 surface defined over a number field $K$.
Then there exist a positive integer $n$ and a finite
extension $L/K$ such that the $L$-rational points of
$S^{[n]}$ are Zariski dense.
\end{theo}

{\em Proof.} By Theorem  \ref{theo:fibrations}
every K3 surface $S$ is 
dominated by an elliptic fibration.
By Remark \ref{rem:main}, we may assume that the class
of the fiber is not big. Therefore it has self-intersection
zero which implies that $S$ is an elliptic K3 surface. 
In this case, the main theorem of \cite{bogomolov-tschi-2}
implies the theorem  with $n=1$. $\square$

\begin{exam}
Let $S$ be a K3 surface of degree 2. Then rational points on 
$S^{[2]}$ are potentially dense. 

Indeed, let $g$ be the polarization. 
By Theorem  \ref{theo:fibrations}, there 
exists an irreducible elliptic curve 
$E\subset S$ such that $g-[E]$ is effective. 
If $g=[E]$ the assertion follows from 
Remark \ref{rem:main}.
Otherwise, 
$$
\langle E, E\rangle < \langle g, E\rangle  < \langle g, g\rangle ,
$$
which implies that $\langle E, E\rangle =0$. 
Then the assertion holds with $n=1$ by \cite{bogomolov-tschi-2}.
\end{exam}

\section{Potential density on $S^{[2]}$}
\label{sect:s2}

Given a fixed K3 surface it is a natural problem to
determine the smallest possible $n$ for which the theorem holds. 
(Of course, we expect that we can always take $n=1$!)
As we have seen, the key to proving potential density
is the existence of abelian fibrations on $S^{[n]}$.

The intersection form on the Picard group of $S$
is an integer-valued nondegenerate quadratic form,
denoted $\langle,\rangle$.  
We recall that the Picard group of $S^{[n]}$ is also equipped
with a natural integer-valued nondegenerate quadratic
form $\left(,\right)$, the Beauville form \cite{beauville-83}.
With respect to this form, we have an orthogonal direct
sum decomposition
$$\Pic(S^{[n]})=\Pic(S) \oplus_{\perp} {\bZ}e,$$
where $\left( e,e \right)=-2(n-1)$ and $2e$ is the class
of the diagonal (more precisely, the nonreduced subchemes in $S^{[n]}$.)

On the K3 surface $S$, the Picard group together with the
quadratic form control much of the geometry of $S$. For example, 
if the quadratic form represents zero, then $S$ admits
an elliptic fibration over $\bP^1$. A naive question would be
whether the analog holds for $S^{[n]}$ with $n\ge 2$.
More precisely, if the Beauville form represents zero, is
$S^{[n]}$ birational to an abelian fibration over $\bP^n$ 
(see \cite{hassett-tschi})?
Note that the Beauville form of $S^{[2]}$ represents zero 
if and only if the intersection form on $\Pic(S)$ represents
$2m^2$ for some $m\in \bZ$.

\begin{prop}\label{prop:2n2}
Let $S$ be a generic K3 surface of degree $2m^2$ with
$m>1$. 
Then $S^{[2]}$ is isomorphic to an abelian surface
fibration over $\bP^2$. 
\end{prop}

{\em Proof.}
We first consider the case $m=2$. 
We asssume that the polarization on $S$ is very ample
and that $S$ does not contain a line.
Then $S$ can be represented as a complete intersection of
a three-dimensional space ${\cal I}_{S}(2)$ of quadrics in $\bP^5$. 
An element of $S^{[2]}$ spans a line $\ell\in \bP^5$
and a two dimensional subspace  of ${\cal I}_{S}(2)$   
contains $\ell$.
In this way, we obtain a morphism 
$$
a\,:\, S^{[2]}\ra \bP^2\simeq \bP({\cal I}_S(2)^*).
$$
The generic fiber of $a$ is an abelian surface; the variety of
lines on a smooth complete intersection of two quadrics in $\bP^5$
is a principally polarized abelian surface (see \cite{GH}, p. 779).  
Notice that $a$ is induced by the sections of $f_8-2e$, 
where $f_8$ is the polarization of degree 8.

When $m>2$ the proof consists of three steps:
\begin{enumerate}
\item{construct special K3 surfaces $S$ so that
$S^{[2]}$ admits a natural involution;}
\item{show directly that some of these special K3 surfaces admit an abelian
surface fibration and a polarization of degree $2m^2$;}
\item{verify that 
this abelian surface fibration deforms to the 
Hilbert scheme of a generic K3 surface of degree $2m^2$.}
\end{enumerate}

We begin with a construction of Beauville and Debarre \cite{debarre}.
Let $S\subset \bP^3$ be a smooth quartic hypersurface;  
in particular, $S$ is a K3 surface and the corresponding polarization
is denoted $f_4$.  Then there is a birational involution
$$
j\,:\, S^{[2]} \dashrightarrow S^{[2]}
$$
defined on an open subset of $S^{[2]}$ by the rule $j(p_1+p_2)=p_3+p_4$,
where $p_1,p_2,p_3,$ and $p_4$ are collinear points on $S$.  
This is a morphism provided that $S$ does not contain a line.  
The action of $j$ on the Picard group of $S^{[2]}$ is given by
$$
j^*x=-x+\left( f_4-e,x \right) (f_4-e).
$$

Next, we consider some special quartic K3 surfaces.  Let $S$ be a K3 surface
with Picard group generated by the ample class $f_4$ and a second
class $f_8$ satisfying
$$
\begin{array}{c|rr}
 & f_4 & f_8 \\
  \hline
  f_4 & 4 & k \\
  f_8 & k & 8 
  \end{array}
$$
where $k>7$.  Such K3 surfaces are parametrized by a nonempty analytic
open subset of an irreducible variety of dimension 18.  This follows
from the Torelli theorem, surjectivity of Torelli,
and the structure of the cohomology lattice of K3 surfaces (see
\cite{looijenga-peters} Theorem 2.4
and \cite{beauville-85}).  Note that $f_4$ is very ample and that 
the image is a smooth quartic surface not 
containing a line \cite{saint-donat};
here we are using the fact that $k\ne 6$.  
Furthermore, the same reasoning shows that $f_8$ is very ample and
the image does not contain a line, provided that $f_8$
is ample.  (Here we are using the fact that $k\ne 7$.)
If $f_8$ were not ample then
$\langle f_8,C\rangle\le 0$ for some $(-2)$-curve $C$
(see \cite{looijenga-peters} 1.6).  
Clearly $\langle f_8,C\rangle\ne
0$ and if
$\langle f_8,C\rangle <0$ then
the Picard-Lefschetz reflection
$\rho(f_8)=f_8+\langle f_8,C\rangle C$ and $f_4$ generate a sublattice with
discriminant greater than $32-k^2$, which is impossible.  
Our argument in the $m=2$ case shows that the $S^{[2]}$ admits
an abelian surface fibration, induced by the line bundle $f_8-2e$.  
Composing with the involution $j$, we obtain a second elliptic
fibration, induced by 
$$j^*(f_8-2e)=2e-f_8+\left(f_8-2e,f_4-e\right) (f_4-e)=(k-4)f_4-f_8-(k-6)e.$$
Let $g=(k-4)f_4-f_8$ and $m=k-6$ so that $\langle g,g\rangle =2(k-6)^2=2m^2$
and $j^*(f_8-2e)=g-me$.  
Note that $g$ is effective on $S$. 

We turn to the last step.  Let ${\cal S} \ra \Delta$ be a general
deformation of $S$ for which $g$ remains algebraic.  The class
$g$ restricts to a polarization on the generic fiber, since it
has Picard group of rank one.  The class $g-me$ is algebraic (and nef)
on the generic fiber of ${\cal S}^{[2]} \ra \Delta$.  
Using deformation theory (see \cite{hassett-tschi} and \cite{ran} Cor. 3.4),
we find that the generic fiber also admits an abelian fibration
with base ${\bP^2}$, induced by the sections of the line bundle
$g-me$.  We are using the fact that the 
abelian surface fibration on $S^{[2]}$ is {\em Lagrangian};  see
\cite{hassett-tschi} for the fourfold case and \cite{matsushita2}
more generally.
$\square$

\begin{rem}
Unfortunately, our argument gives little information about how the
abelian fibration degenerates for nongeneric K3 surfaces of degree
$2m^2$ with $m>2$.  
A more precise description would follow from the conjectures
of \cite{hassett-tschi}.  
\end{rem}

\begin{theo}
\label{theo:s2}
Let $S_8$ be a K3 surface of degree 8, defined
over a number field $K$, embedded in projective space
$\bP^5$ as a complete intersection of 3 quadrics and not containing
a line. Then rational points on $S_8^{[2]}$ are potentially dense. 
The same result holds for a generic K3 surface of degree $2m^2$.
\end{theo}

{\em Proof.}
We apply Theorem \ref{theo:dense}, using the first set of assumptions.
We use the abelian fibrations constructed in Proposition \ref{prop:2n2}.

Let $g$ be the homology class of an irreducible elliptic curve
(see Theorem  \ref{theo:fibrations}). 
We verify that $g*g$ intersects the class of a fiber
positively.

We need to compute the intersection on $S^{[2]}$ of 
$(f-me)\cdot(f-me)\cdot (g*g)$, where  $f$ and $g$ are divisor classes
on $S$. Let $\Sigma$ be the class of subschemes 
containing a fixed point $p\in S$;
note that these subschemes are parametrized by 
the blow-up of $S$ at $p$. In particular, 
$(f-me)\cdot(f-me)\cdot\Sigma=\langle f,f\rangle -m^2$ 
(because $e$ restricts
to the exceptional divisor of the blown-up K3 surface).
We also have 
\begin{eqnarray*}
g*g  & = & g\cdot g - \langle g,g\rangle \Sigma,\\
f\cdot f\cdot g\cdot g &  = & \langle f,f\rangle \langle g,g\rangle 
+2\langle f,g\rangle ^2,\\
f\cdot e\cdot g\cdot g & =  & 0,\\
e\cdot e\cdot g\cdot g & =  & -2\langle g,g\rangle .
\end{eqnarray*}
Finally, 
we obtain 
$$
(f-me)\cdot(f-me)\cdot (g*g)= 2\langle f,g\rangle^2-
m^2\langle g,g\rangle .
$$
In our case, $f=f_{2m^2}$, $g$ is the class of the elliptic
curve. 
To verify the hypothesis of the Theorem \ref{theo:dense}, 
we need $2\langle f_{2m^2},g\rangle^2 > 
m^2 \langle g,g\rangle $. Since $\langle g,g\rangle >0$ 
we are done by the
Hodge index theorem,  
which implies that the determinant of the matrix
$$
\left(
\begin{array}{cc} 2m^2 & 
\langle f_{2m^2},g\rangle  \\ \langle f_{2m^2},g \rangle & 
\langle g,g\rangle \end{array} 
\right)
$$
is negative. 
$\square$

\end{document}